\documentclass{article}
\usepackage{amsmath,  amsthm, amssymb}
\numberwithin{equation}{section}

\begin{document}

\author{Ajai Choudhry}
\title{Quadruples of  integers whose elementary \\symmetric functions are perfect squares}

\date{}
\maketitle
\abstract{While Euler had considered the problem of finding three integers whose elementary symmetric functions are all squares, an extension of the  problem to  four integers has not been discussed  in the mathematical literature till now. In this paper, we describe a method of generating infinitely many quadruples of positive integers whose elementary symmetric functions are all squares}. 

\medskip

\noindent Keywords: four integers whose elementary symmetric functions are all squares.
\medskip

\bigskip

\section{Introduction}\label{intro}
In 1760 Euler considered the problem of finding three distinct positive integers whose elementary symmetric functions are all squares. He found this to be a very difficult problem, and his first method of solving led to very  large integers, with the smallest numerical example  consisting of integers with 13 digits \cite{LE1}. Later on, Euler discovered a simpler method of solving the problem by imposing a simplifying condition \cite{LE2}. Both methods of Euler yield parametric solutions of the problem in terms of polyniomials of high degrees. While Euler did not peform the cumbersome computations, the simplest parametric solutions obtained by his two methods are of degrees 132 and 16, respectively. Recently, Choudhry \cite{Cho} found much simpler parametric solutions of the problem and obtained several numerically  smaller solutions.

In this paper we consider an extension of the problem to four integers, that is, we are concerned with  finding quadruples of  integers whose four elementary symmetric functions are all squares. It seems that till now this problem has not  been considered at all in the mathematical literature. We obtain several numerical examples of such quadruples and also  describe a method of generating infinitely many quadruples of  positive integers whose  elementary symmetric functions are all squares. The  smallest such quadruple  with this property is  $\{ 420, 525, 1344, 1680\}.$

\section{Quadruples of  integers whose  elementary symmetric functions are all squares}\label{quadruples}
If $a_1, a_2, a_3, a_4$ are four integers whose  elementary symmetric functions are all squares and $k$ is a nonzero integer, it is readily seen that  
the four integers $k^2a_1, k^2a_2, k^2a_3, k^2a_4$  also  have the same property. Similarly, if the greatest common divisor of the integers $a_1, a_2, a_3, a_4$ has a squared factor $k^2$, then we may divide each of the integers $a_i, i=1,\ldots, 4$ by $k^2$ to get a new set of four numerically smaller  integers whose  elementary symmetric functions are all squares. A quadruple of integers $a_1, a_2, a_3, a_4$ will be said to  be in a reduced form if the integers $a_i, i=1, \ldots, 4$, do not have a common squared factor. All numerical solutions of the problem  obtained in this paper will be presented in  reduced form. While the integers $a_i, i=1, \ldots, 4$, of a quadruple in reduced form do not have a common squared factor, they may, of course, have a common factor that is not a perfect square. Indeed, in all the numerical solutions that we obtain, the four integers of the quadruple do have a non-squared common factor  greater than 1.

Let $a_1, a_2, a_3, a_4$ be four integers having a common factor $m$ so that $a_i=mx_i$, where $x_i, i=1, \ldots, 4$, are certain integers. The elementary symmetric functions $s_i, i=1, \ldots, 4$, of the integers $a_1, a_2, a_3, a_4$ may now be written as follows:
\begin{align}
s_1 & = m(x_1+x_2+x_3+x_4), \label{sym1} \\
s_2 & = m^2(x_1x_2+x_1x_3+x_1x_4+x_2x_3+x_2x_4+x_3x_4), \label{sym2} \\
s_3 & =m^3(x_1x_2x_3+x_1x_2x_4+x_1x_3x_4+x_2x_3x_4), \label{sym3} \\
s_4 & = m^4x_1x_2x_3x_4. \label{sym4} 
\end{align}
Since $s_3 = m^3\{x_1x_2(x_3+x_4)+x_3x_4(x_1+x_2)\}$, if we choose integers $x_i, i=1, 2, 3, 4$, such that the auxiliary conditions
\begin{equation}
x_1x_2 =x_3x_4 =t^2, \label{aux1}
\end{equation}
are satisfied, then $s_4=m^4t^4$, and we get $s_3=m^3t^2(x_1+x_2+x_3+x_4)$, so that both $s_1$ and $s_3$ will become squares if we take $m=x_1+x_2+x_3+x_4$. 

A solution of Eqs. \eqref{aux1} is readily obtained and may be written as follows:
\begin{equation}
x_1=fg^2, \quad x_2=fh^2, \quad x_3=gf^2, \quad x_4=gh^2, \quad t=fgh, \label{valx}
\end{equation}
where $f, g$, and $h$ are arbitrary integer parameters. Further, using the values of $x_i$ given by \eqref{valx}, we get
\begin{equation}
m = x_1+x_2+x_3+x_4 = (f + g)(fg + h^2). \label{valm} 
\end{equation}

With the values of $x_1,\ldots, x_4$ and $m$ given by \eqref{valx} and \eqref{valm}, respectively, the three elementary functions $s_1, s_3$, and $s_4$ become perfect squares. In fact, we get,
\begin{equation}
\begin{aligned}
s_1 & = (f + g)^2(fg + h^2)^2, \\
s_3 & = f^2g^2h^2(f + g)^4(fg + h^2)^4, \\
s_4 & = f^4g^4h^4(f + g)^4(fg + h^2)^4.
\end{aligned}
\end{equation}
Further,  on substituting the values of $x_1,\ldots, x_4$ and $m$ in \eqref{sym2}, we get 
\begin{equation}
s_2 = fg(f + g)^2(fg + h^2)^2(h^4 + (f + g)^2h^2 + f^2g^2). \label{vals2}
\end{equation}

 It follows from Eq. \eqref{vals2} that $s_2$ will also become a perfect square if we find  integers $f, g, h$, and $y$ satisfying the condition
\begin{equation}
y^2=fgh^4 + fg(f + g)^2h^2 + f^3g^3. \label{condy}
\end{equation}
Each solution of Eq. \eqref{condy} may be used to generate a quadruple of  integers $a_i, i=1, \ldots, 4$, whose elementary symmetric functions are all perfect squares. 

In the next subsection, we present  numerical solutions of Eq. \eqref{condy} and the quadruples thus generated, and in subsection \ref{quadruplesec} we give a method of generating infinitely many solutions of Eq. \eqref{condy} that will yield infinitely many quadruples of positive integers whose elementary symmetric functions are all squares.

\subsection{Numerical examples of quadruples} \label{quadruplesnumeric} Computer trials over the range $f+g+h \leq 400$ yielded a large number of  integer solutions of Eq. \eqref{condy} and each such solution generated a quadruple of integers whose elementary symmetric  functions are all squares.  When we expressed these quadruples in reduced form, we found only five distinct quadruples in which all the four integers were less than $10^7$. In Table \ref{quadruplesa1234} we present these  five quadruples in reduced form with the values of the integers written in ascending order as well as the values of the triple $(f, g, h)$ that yielded each quadruple. We  also note that each  quadruple in reduced form was generated by several values of the triple $(f, g, h)$. In Table \ref{quadruplesa1234} we have presented the values of the triple $(f, g, h)$ in which the sum $f+g+h$ has the smallest value. 
\begin{table}
\begin{center}
\begin{tabular}{|c|c|} \hline
$(f, g, h) $  & $ a_1, a_2, a_3, a_4$ \\
\hline
(4, 5, 8)     & 420, 525, 1344, 1680 \\
(5, 8, 15)    & 27560, 44096, 155025, 248040 \\
(1, 40, 3)    & 369, 1640, 14760, 65600 \\
(16, 36, 21)  & 287924, 376064, 647829, 846144 \\
(8, 13, 65)   & 49728, 80808, 2020200, 3282825 \\
\hline
\end{tabular}
\caption{Quadruples of integers $a_1, a_2, a_3$, and $a_4$} {whose  elementary symmetric functions are all squares.}
\label{quadruplesa1234}
\end{center}
\end{table}

The sum of the four integers of the smallest quadruple in Table \ref{quadruplesa1234} namely, \{420, 525, 1344, 1680\} is $3969=63^2$. We accordingly performed a further computer search for all quadruples $a_1, a_2, a_3$, and $a_4$ of positive integers whose sum is $m^2$ with $m \leq 63$ and whose elementary symmetric functions are all squares. No smaller quadruple was found. Thus, the smallest such quadruple is  \{420, 525, 1344, 1680\} as already mentioned in the introduction. 

Finally, we note that all the four integers of each quadruple $a_1, a_2, a_3, a_4$ listed in Table \ref{quadruplesa1234} kave a common factor greater than 1. This is also true for all other quadruples that we found by the computer trials mentioned above and which are not listed in Table \ref{quadruplesa1234}.

\subsection{Method of generating infinitely many quadruples}\label{quadruplesec}
We now describe a method of generating infinitly many quadruples of positive integers whose elementary symmetric functions are all squares. Since the first entry in Table \ref{quadruplesa1234}  has $(f, g, h)=(4, 5, 8)$, we now take $(g, h)=(5, 8)$ when Eq. \eqref{condy} reduces to the following cubic equation:
\begin{equation}
y^2 = 445f^3 + 3200f^2 + 28480f. \label{ecfy}
\end{equation}
Eq. \eqref{ecfy} represents an elliptic curve and there is a rational point on it given by $(f, y)=(4, 440)$. In fact, the birational transformation defined by
\begin{equation}
f = (4X)/445 - 12/5,\quad  y = (8Y)/445, \label{biratfy}
\end{equation}
and 
\begin{equation}
X = (445f)/4 + 267, \quad  Y = (445y)/8, \label{biratXY}
\end{equation}
transforms the cubic curve \eqref{ecfy} to the standard Weierstrass form
\begin{equation}
Y^2 = X^3 - X^2 + 578767X - 173493663. \label{ecXY}
\end{equation}

Corresponding to the known rational point (4, 440) on the curve \eqref{ecfy}, there is a rational point $P$ on the elliptic curve \eqref{ecXY} whose coordinates are $(712, 24475)$. We can now find more rational points on the elliptic curve \eqref{ecXY} using the group law. The point $2P$ on the curve \eqref{ecXY} is $(200517/484, -124731987/10648).$ Since the point  $2P$ has rational coordinates, it follows from the Nagell-Lutz theorem \cite[p. 56]{ST} on elliptic curves that $P$ is not a point of finite order. 

The software SAGE determined that the rank of the elliptic curve \eqref{ecXY} is 2, with its two generators being $(712, 24475)$ and $(1909/4 , 116261/8)$. This further confirms the existence of infinitely many rational points on the curve \eqref{ecXY}. We can thus compute infinitely many rational points on the elliptic curve \eqref{ecXY} using the group law, and corresponding to these points, we can find infinitely many rational points on the cubic curve \eqref{ecfy}. 

We note here that the right-hand side of Eq. \eqref{ecfy} may be written as $f(445f^2 + 3200f + 28480)$. The discriminant of the quadratic factor is negative, and hence the quadratic factor is necessarly positive for rational values of $f$. For any nontrivial rational solution of \eqref{ecfy}, both sides are positive, and hence the factor $f$ on the right-hand side must be necessarily positive. Since we had fixed $(g, h)=(5, 8)$, the method described above ensures that we  get infinitely many rational solutions of Eq. \eqref{condy} with all three $f, g$, and $h$ being positive. .

Since the right-hand side of Eq. \eqref{condy} is a homogeneous function of degree 6 in the parameter $f, g$,  and $h$, we can multiply each of the rational numbers $f, g$, and $h$ by a suitable positive integer to get infinitely many  solutions in positive integers of Eq. \eqref{condy}. Now, using the relations \eqref{valx} and \eqref{valm}, and then the relations $a_i=mx_i, i=1,\ldots, 4$, we can generate infinitely many quadruples of positive integers $\{a_1, a_2, a_3, a_4\}$ whose elementary symmetric functions are all squares. 

As a numerical example, the point $2P$ on the curve \eqref{ecXY} yields the point $(801/605,$ $ -1401483/6655)$ on the curve \eqref{ecfy}. We thus get a solution of Eq. \eqref{condy} with $(f, g, h)=(801/605, 5, 8)$. On multiplying each of these numbers by 605, we get  $(f, g, h)=(801,$  $3025, 4840)$. These values of $f, g$, and $h$ yield a quadruple of integers which, in reduced form, may be written as follows:
\[
\{20975970565170,\quad  79216368239250, \quad 202793902692480, \quad   765857123152000\}.
\]

In Section \ref{quadruplesnumeric} we had noted that the same reduced quadruple may be generated by multiple values of $(f, g, h)$.  Accordingly, we  now show that the above procedure actually generates infinitely many reduced quadruples. Each reduced quadruple determines three ratios $a_1/a_2, a_1/a_3$, and $a_1/a_4$. Two quadruples will be identical only if these three ratios are the same for both of them though these three ratios  may  appear in a different order for the two quadruples. Using the relations $a_i=mx_i, i=1,\ldots, 4$, and \eqref{valx}, and noting that we have already fixed $(g, h)=(5, 8)$, these three ratios are $\{25/64, 5/f, 5f/64\}$. It is now readily 
seen that if a certain triple of ratios is generated by $f=f_0$, it can only be generated by one more value of $f$, namely $f=64/f_0.$ Thus, any reduced quadruple can at most be generated twice by the above procedure. Since we can obtain infinitely many rational values of $f$, we are assured of infinitely many distinct reduced quadruples whose elementary symmetric functions are all squares.

We now note that in Eq. \eqref{condy}, we could also fix the values of $f$ and $h$ instead of $g$ and $h$. Thus, keeping in view the values of $(f, g, h)=(4, 5, 8)$ that generate the first solution listed in Table \ref{quadruplesa1234}, we may fix $(f, h) =(4, 8)$,  when Eq. \eqref{condy} reduces to the cubic curve
\begin{equation}
y^2 = 320g^3 + 2048g^2 + 20480g. \label{ecgy}
\end{equation}.
We have verified that Eq. \eqref{ecgy} represents an elliptic curve of rank 1, and we can proceed as above, to generate new examples of quadruples of integers whose elementary symmetric functions are all squares. Similarly, on taking $(f, g)=(4, 5)$, Eq. \eqref{condy} reduces to the following quartic model of an elliptic curve:
\begin{equation}
y^2 = 20h^4 + 1620h^2 + 8000. \label{echy}
\end{equation}
The rank of the elliptic curve  \eqref{echy} is also 1, and it can also be used to generate infinitely many quadruples whose elementary symmetric functions are all squares. 

Finally, we can follow the same procedure using any of the other numerical solutions listed in Table \ref{quadruplesa1234} to generate more examples of quadruples with the desired property.

\section{Some open problems} Since we can obtain several elliptic curves of positive rank that yield the desired quadruples, it is possible that there exists a parametric solution of the problem. Moreover, all the quadruples in reduced form that we found have a common factor greater than 1. It will be therefore be interesting to find a parametric solution, and also investigate the existence of numerical quadruples in which the four integers have no common factor. 

Finally, there remains the open problem of finding sets of more than 4 integers such that all the elementary functions of such sets of integers are perfect squares. We leave these problems for future investigations.

\noindent Ajai Choudhry, 13/4 A Clay Square, Lucknow - 226001, India.

\smallskip

\noindent E-mail address: ajaic203@yahoo.com

\end{document}